\documentclass[12pt]{article} 
\usepackage{amstext,amsfonts,graphics}
\newtheorem{thm}{Theorem} \newtheorem{lemma}[thm]{Lemma}
 \newtheorem{cor}[thm]{Corollary}
\newenvironment{pf} {\noindent{\sc Proof. }}{{\hfill
$\Box$}\par\vskip2\parsep} 
\newenvironment{pfof}[1]
{\par\vskip2\parsep\noindent{\sc Proof of\ #1. }}{{\hfill $\Box$}
  \par\vskip2\parsep}

\renewcommand{\P}{{\bf P}}
\newcommand{\Z}{{\mathbb Z}}
\newcommand{\py}{{\bf P}}
\newcommand{\Q}{{\bf Q}}

\newcommand{\E}{{\bf E}}
\newcommand{\n}{N}

\newcommand{\dof}{\bf\boldmath} 
\newcommand{\ind}{\mbox{\boldmath $1$}}

\newcommand{\sig}{\sigma}
\newcounter{mycount}
\newenvironment{mylist}{\begin{list}{(\roman{mycount})}%
{\usecounter{mycount}\itemsep 0pt}}{\end{list}}

\title{An invariant of finitary codes with finite expected
 square root coding length}
\date{July 7, 2003}

\author{Nate Harvey\thanks{Department of Mathematics, UC Berkeley, CA
    94720, USA.} \and Yuval Peres\thanks{Departments of Statistics and
    Mathematics, UC Berkeley, CA 94720, USA. {\tt
      peres@stat.berkeley.edu}. Research supported by NSF Grants
DMS-0104073, CCR-0121555 and a Miller Professorship at UC Berkeley.}}
\begin{document}

\maketitle
\begin{abstract}
  Let $p$ and $q$ be probability vectors with the same entropy $h$.
  Denote by $B(p)$ the Bernoulli shift indexed by $\Z$
  with marginal distribution $p$. Suppose
  that $\varphi$ is a measure preserving homomorphism from $B(p)$ to
  $B(q)$.  We prove that if the coding length of
  $\varphi$ has a finite $1/2$ moment, then $\sigma_p^2=\sigma_q^2$, 
   where $\sigma_p^2=\sum_i p_i(-\log p_i-h)^2$ is
  the {\em informational variance\/} of $p$. In this result,
 the $1/2$ moment cannot be replaced by a lower moment. On the other
  hand, for any $\theta<1$, we exhibit probability vectors $p$ and $q$
  that are not permutations of each other, such that there exists a
  finitary isomorphism $\Phi$ from $B(p)$ to $B(q)$ where the coding
  lengths of $\Phi$ and of its inverse have a finite $\theta$ moment.
  We also present an extension to ergodic Markov chains.
\end{abstract}

\section{Introduction}
Let ${\bf A}= \{\alpha_0, \ldots ,\alpha_{a-1}\}$ be a finite alphabet
and $p=(p_0, \ldots ,p_{a-1})$ a probability vector with 
entropy $h(p) = \sum_{i=0}^{a-1} -p_i \log(p_i)$.
Consider the Bernoulli shift
 $B(p)=(X,{\cal A},\P , T)$, where $X = {\bf A}^{\Z}$
is equipped with the product $\sigma$-algebra ${\cal A}$,
the product measure $\P = p^{\mathbb Z}$ and the left shift $T$. 
Let ${\bf B}=\{\beta_0, \ldots ,\beta_{b-1}\}$ be another finite alphabet,
and $q=(q_0, \ldots ,q_{b-1})$ a probability vector; denote by 
$B(q)=(Y,{\cal B}, \Q, T)$ the corresponding  Bernoulli shift.
A {\dof homomorphism} $\varphi$ from $B(p)$ to $B(q)$ is a
measurable map from $X$ to $Y$, defined $\P$-a.e., such that 
$\P\varphi^{-1}=  \Q$ and  $\varphi T=T \varphi \, $ $ \, \P$-a.e.. 
An {\dof isomorphism} is an invertible homomorphism.
A homomorphism $\varphi$ from $B(p)$ to $B(q)$ is
{\dof finitary} if there exists a set $W \subseteq X$ with
$\P(W)=1$,  that has the following property: 
 for all $x \in W$ there exists $n=n(x)$ such that if
$\widetilde{x} \in W$ and 
$\widetilde{x}_i=x_i$ for all $ -n \le i \le n$,   then
$(\varphi(x))_{0}=(\varphi(\widetilde{x}))_{0}$.  We write
$\n_{\varphi}(x)$ for the minimal such $n$, and call $\n_{\varphi}(x)$
the {\dof coding length} of $\varphi$.  A {\dof finitary isomorphism} is an
invertible finitary homomorphism whose inverse is also finitary.

 By the Kolmogorov-Sinai Theorem (see, e.g., \cite{pete}),
if $B(p)$ and $B(q)$ are isomorphic, then $h(p)=h(q)$. 
The converse was established by Ornstein~\cite{o}.
Keane and Smorodinsky~\cite{ks:iso} proved that if $h(p)=h(q)$, then
there exists a finitary isomorphism from $B(p)$ to $B(q)$.
Parry \cite{parry:II} and Schmidt \cite{schmidt:one}  showed
 that if a  finitary isomorphism from $B(p)$ to $B(q)$
has finite expected coding length in both directions,
then $p$ and $q$ must be permutations of each other.

In this paper, we prove that  
the {\dof informational variance} of $p$,
$$
\sigma^2_p=\sum_{i=0}^{a-1}p_i\Bigl(-\log (p_i)-h(p)\Bigr)^2 \,
$$
is an invariant of isomorphisms $\varphi$ that satisfy
$\E\Bigl(\n_\varphi^{1/2}\Bigr)< \infty$.  
More precisely:

\begin{thm}
\label{strong}
Let $p$ and $q$ be probability vectors that satisfy
$h(p)=h(q)$ and $\sigma^2_p \neq \sigma^2_q$. Then there
exists a constant $c_{p,q}>0$ such that for any
finitary homomorphism  $\varphi$ from $B(p)$ to $B(q)$,
we have 
$$ 
 \liminf_{n \rightarrow \infty} 
 \frac {\E ( \n_\varphi \wedge n) }{\sqrt{n}}  \ge c_{p,q} 
$$
and consequently, $\E\Bigl(\n_\varphi^{1/2}\Bigr)= \infty$.  
\end{thm}
(Here and throughout, $\E$ denotes expectation with respect to
$\P=p^\Z$.)

The exponent $1/2$ in the theorem is sharp, since Meshalkin \cite{mesh}
(see \S 3) constructed a finitary isomorphism $\varphi$ from $B(p)$
for $p=\Bigl({1 \over 2},{1 \over 8},{1 \over 8},{1 \over 8}, {1
    \over 8}\Bigr)$ to $B(q)$ for 
$q=\Bigl({1 \over 4},{1 \over 4},{1 \over 4},{1 \over 4}\Bigr)$, 
where $\P[\n_\varphi>k]$ equals the probability
that a simple random walk remains positive for $k$ steps.
Thus for Meshalkin's code, $0< \lim_k \P[\n_\varphi>k] \sqrt{k} <\infty$,
whence $\E\Bigl(\n_\varphi^{\theta}\Bigr)< \infty$
for all $\theta<1/2$. Clearly $\sigma_q<\sigma_p$ in this case,
so Meshalkin's code is essentially optimal.  

The assumption that $\sigma^2_p \neq \sigma^2_q$ in Theorem
\ref{strong} cannot be dropped, as shown by our next result.
\begin{thm}
\label{one}
For any $0<\theta<1$, there are probability vectors $p$ and $q$ 
where $p$ is not a permutation of $q$, such that
there exists a finitary isomorphism $\Phi$ from $B(p)$
to $B(q)$ that satisfies $\E(\n_\Phi^\theta) < \infty$ and 
$\E_\Q(\n_{\Phi^{-1}}^\theta) < \infty$.
\end{thm}

Theorem \ref{strong} is proved in the next section.  In \S 3 we 
recall Meshalkin's isomorphism, and describe an
adaptation of Meshalkin's code which motivates Theorem \ref{one}.  
In \S 4 we define a class of matchings useful for
the proof of Theorem \ref{one}, and in \S 5 we prove the theorem.  In
\S 6 we define informational variance for ergodic Markov chains, and
present an extension of Theorem \ref{strong} to this setting.

\section{Proof of Theorem 1}
With the notation of the introduction in force, we may assume that
the probability vectors $p$ and $q$ satisfy
$p_i>0$ for all $0 \le i <a$ and $q_j>0$ for all $0 \le j <b$.
 Let $\varphi$ be a finitary homomorphism  from $B(p)$ to $B(q)$.
For $x =(x_k)_{k \in \Z}\in X$, write  
$X_i(x) =-\log(p(x_i))-h(p)$, where $p(\alpha_j)=p_j$
for any $j$.  Similarly, if 
$ \varphi(x)=y=(y_k)_{k \in \Z}\in Y$, let $Y_i(x)=
-\log(q(y_i))-h(q)$.  Since $\varphi\P^{-1}=\Q$, it
follows that $\E(X_i)=\E(Y_i)=0$.
Let $S_n = \sum_{i=1}^n X_i$ and $R_n =  \sum_{i=1}^n Y_i$. 
Write $t^+=\max\{t,0\}$.
\begin{lemma}
\label{var}
If $\sigma^2_p \neq \sigma^2_q$, then 
$$
\liminf_{n \to \infty} \frac{1}{\sqrt n}\E (R_n-S_n)^+ \ge
\frac{|\sig_q-\sig_p|}{\sqrt{2\pi}} \, .
$$
\end{lemma}
\begin{pf}
By a version of the central limit theorem
(see \cite{st}, Cor.\ 2.1.9),
$$
\lim_{n \to \infty} \E\Bigl(\frac{R_n^+}{\sqrt{n}}\Bigr) =
\frac{\sig_q}{\sqrt{2\pi}}\int_0^\infty t e^{-\frac{t^2}{2}} \, dt 
=\frac{\sig_q}{\sqrt{2\pi}} \,,
$$
and similarly
$$
\lim_{n \to \infty} \E\Bigl(\frac{S_n^+}{\sqrt{n}}\Bigr) 
=\frac{\sig_p}{\sqrt{2\pi}} \,.
$$

Since $  (R_n-S_n)^+ \ge R_n^+-S_n^+$, we infer that
\begin{equation} \label{slarge}
\liminf_{n \to \infty} \frac{1}{\sqrt n}\E (R_n-S_n)^+ \ge
\frac{\sig_q-\sig_p}{\sqrt{2\pi}} \, 
\end{equation}
and similarly
\begin{equation} \label{rlarge}
\liminf_{n \to \infty} \frac{1}{\sqrt n}\E (S_n-R_n)^+ \ge
\frac{\sig_p-\sig_q}{\sqrt{2\pi}} \, .
\end{equation}

If $\sig_q>\sig_p$, then (\ref{slarge}) proves the lemma.
In the remaining case,  $\sig_p>\sig_q$, 
the assertion of the lemma follows from (\ref{rlarge})
by taking expectations in the identity
$$
(R_n-S_n)^+ =(R_n-S_n) + (S_n-R_n)^+ \,.
$$
\end{pf}

\begin{lemma}
\label{short}
 Let $\varphi$ be a finitary homomorphism  from $B(p)$ to $B(q)$.
Denote $\lambda_q = \max \{-\log(q_j) \,: \, 0 \leq j \leq b-1\}$.
Then for all $n$,
$$
 \E (R_n-S_n)^+ \le  2\lambda_q \E(\n_{\varphi} \wedge n) \,.
$$
\end{lemma}
\begin{pf}
Let
$$
I_n=I_n(x)= \Bigl\{i \in \{1,\ldots, n\}:\n_{\varphi}(T^i x) >
\min\{i,n+1-i\} \Bigr\}
$$
and denote $J_n=\{1,\ldots,n\} \setminus I_n$. Observe that
\begin{equation} \label{ej}
\E|I_n| = \sum_{i=1}^n \py(i \in I_n)
\le 2\sum_{i=1}^n
\py(\n_\varphi \geq i)  \le 2\E(\n_{\varphi} \wedge n) \,.
\end{equation}

Fix $x \in X$ and let $y=\varphi(x)$.
Since 
$$
\Bigl\{\widetilde{x} \in X : 
(\widetilde{x}_1, \ldots, \widetilde{x}_n)  =  (x_1, \ldots, x_n) \Bigr\} 
  \subset 
\varphi^{-1} \Bigl\{\widetilde{y} \in Y :  
\widetilde{y}_j=y_j  \; \forall j \in J_n \Bigr\} \, ,
$$
it follows that
$$
\P \Bigl\{\widetilde{x} \in X : 
(\widetilde{x}_1, \ldots, \widetilde{x}_n)  =  (x_1, \ldots, x_n) \Bigr\} 
  \le \Q \Bigl\{\widetilde{y} \in Y :  
\widetilde{y}_j=y_j  \; \forall j \in J_n \Bigr\} \, .
$$

Taking logarithms, this implies that 
$$
\sum_{k=1}^n \log p(x_k) \le
\sum_{k=1}^n \log q(y_k) -
\sum_{i  \in I_n } \log q(y_i) \le \sum_{k=1}^n \log q(y_k) 
 + \lambda_q |I_n| \,.
$$

Since $h(p)=h(q)$, we deduce from the last equation and the definitions of
$R_n$ and $S_n$ that 
$R_n-S_n \leq \lambda_q |I_n|$, whence by (\ref{ej}),
$$
 \E (R_n-S_n)^+ \le  \lambda_q \E |I_n|
\le 2\lambda_q \E(\n_{\varphi} \wedge n) \,.
$$

\end{pf}

\begin{pfof}{Theorem \ref{strong}}
Lemmas \ref{var} and \ref{short} imply that
\begin{equation} \label{conc1} 
 \liminf_{n \rightarrow \infty} 
 \frac {\E ( \n_\varphi \wedge n) }{\sqrt{n}}  \ge 
\frac{|\sig_q-\sig_p|}{2\lambda_q \sqrt{2\pi}} >0 \, ,
\end{equation}
so it only remains to verify the final assertion of the theorem.

Observe that $\n_\varphi \wedge n \le  \sqrt{\n_{\varphi}n} \;$
and $\: (\n_\varphi \wedge n)/\sqrt{n} \to 0 \; \; \P$-a.e.

\noindent{} If we had
$\E\Bigl(\sqrt{\n_{\varphi}}\Bigr)< \infty$, then we could deduce 
by dominated convergence that
$\E( \n_\varphi \wedge n)/\sqrt{n} \to 0$, 
which contradicts (\ref{conc1}). Thus 
$\E\Bigl(\sqrt{\n_{\varphi}}\Bigr) =\infty$.
\end{pfof}
A similar idea was used in a different context by Liggett \cite{lig:head}.
\section{Motivating examples and heuristics}

\noindent{\bf Meshalkin's coding}

First, we briefly recall the Meshalkin isomorphism \cite{mesh}.  Let
$B(r)$ be the Bernoulli shift on the alphabet 
${\bf A}_1=\{\alpha_1,\ldots,\alpha_5\}$ for 
$r=\left({1 \over 2}, {1 \over 8},{1 \over 8},{1 \over 8},{1\over 8}\right)$ 
and let $B(s)$ be the
Bernoulli shift on the alphabet ${\bf B}_1=\{\beta_1,\ldots,\beta_4\}$
for $s=\left({1 \over 4},{1 \over 4},{1 \over 4},{1 \over 4}\right)$.
We represent the symbols of ${\bf A}_1$ as

\enspace

\begin{picture}(320,30)
\put(0,20){$\alpha_1$}
\put(15,20){$=$}
\put(28,20){$0$}
\put(38,20){,}
\put(48,20){$\alpha_2$}
\put(63,20){$=$}
\put(76,20){$1$}
\put(76,10){$0$}
\put(76,0){$0$}
\put(86,20){,}
\put(96,20){$\alpha_3$}
\put(111,20){$=$}
\put(124,20){$1$}
\put(124,10){$0$}
\put(124,0){$1$}
\put(134,20){,}
\put(144,20){$\alpha_4$}
\put(159,20){$=$}
\put(172,20){$1$}
\put(172,10){$1$}
\put(172,0){$0$}
\put(182,20){,}
\put(191,20){$\alpha_5$}
\put(206,20){$=$}
\put(219,20){$1$}
\put(219,10){$1$}
\put(219,0){$1$}
\put(229,20){,}
\end{picture}

\enspace

 The symbols of ${\bf B}_1$ are
represented as:

\enspace

\begin{picture}(320,20)
\put(0,10){$\beta_1$}
\put(15,10){$=$}
\put(28,10){$0$}
\put(28,0){$0$}
\put(38,10){,}
\put(48,10){$\beta_2$}
\put(63,10){$=$}
\put(76,10){$0$}
\put(76,0){$1$}
\put(86,10){,}
\put(96,10){$\beta_3$}
\put(111,10){$=$}
\put(124,10){$1$}
\put(124,0){$0$}
\put(134,10){,}
\put(144,10){$\beta_4$}
\put(159,10){$=$}
\put(172,10){$1$}
\put(172,0){$1$}
\put(182,10){,}
\end{picture}

\enspace

The Meshalkin finitary isomorphism $\varphi$ from $B(r)$ to $B(s)$ can be
described in two equivalent ways. Given a sequence
$x=(x_j)_{j \in \Z} \in {\bf A}_1^\Z$, denote by $\ell_i$ the length of the
binary representation of  $x_i \in {\bf A}_1$.  The {\bf
random walk description} of $\varphi$ is obtained by defining,
for each $i$ with $\ell_i=1$,
\begin{equation} \label{rw}
m(i)=\min\Bigl\{ m \ge i : \sum_{j=i}^{m} (\ell_i-2) =0 \Bigr\} \,.
\end{equation}
Observe that $m(\cdot)$ is an injective map from 
$\{i\in \Z : \ell_i=1\}$ onto $\{j\in \Z : \ell_j=3\}$.
For each $i \in \Z$ with $\ell_i=1$, remove the bottom bit from
$x_{m(i)}$ and append it at the bottom of $x_i$. This
produces two symbols from ${\bf B}_1$ that are denoted $y_{m(i)}$
 and $y_i$, respectively. Set $\varphi(x)=y=(y_j)_{j \in \Z}$. 

Alternatively, we have an equivalent {\bf inductive construction} of $\varphi$:

\noindent {\bf Step 1:} For each $i \in {\Z}$ such that 
 $\ell_i=1$ and $\ell_{i+1}=3$, send the bottom bit of $x_{i+1}$
below $x_i$, output the resulting ${\bf B}_1$ symbols and {\dof
remove from consideration} both $i$ and $i+1$.

\enspace
\noindent For each $n \geq 2$, perform:

\noindent {\bf Step $n$:} For all $i \in \Z$ such that 
$\ell_i=1, \ell_{i+n}=3$ and 
$i,i+n$ have not been removed
from consideration, send the bottom bit of $x_{i+n}$ below $x_i$,
output the corresponding ${\bf B}_1$ symbols and {\dof remove from
consideration} both $i$ and $i+n$.

\enspace
\enspace

\noindent {\bf An adaptation of Meshalkin's coding}

Next we describe informally a variant of the coding above, which we will
generalize in \S 5 to prove Theorem \ref{one}.  Consider the random
walk where each increment $X_i$ has $\P(X_i=1)=\P(X_i=3)={1 \over 2}$.
The moment generating function is
\begin{eqnarray*}
\Gamma(z)=\E(z^{X_i})={z+z^3 \over 2}.
\end{eqnarray*}
Consider also the walk where each increment $Y_i$ equals $2$ 
with probability $1$.  This has moment generating function
\begin{eqnarray*}
\Delta(z)=\E(z^{Y_i})=z^2.
\end{eqnarray*}
These walks count the accumulated information for the Bernoulli shifts
$B(r)$ and $B(s)$, where
$r=\left({1 \over 2}, {1 \over 8},{1 \over 8},{1 \over 8},{1 \over 8}\right)$
and $s=\left({1 \over 4},{1 \over 4},{1 \over 4},{1 \over 4}\right)$
as in Meshalkin's coding.
The entropy equality $h(r)=h(s)$ corresponds to the identity
$\Gamma'(1)=\Delta'(1)$ while the inequality of informational variance
corresponds to the inequality $\Gamma''(1) \ne \Delta''(1)$.
The identity 
\begin{eqnarray*}
\Gamma^2(z)-\Delta^2(z)={1 \over 2}\bigg(\Gamma\bigg(z^2\bigg)-\Delta\bigg(z^2\bigg)\bigg)
\end{eqnarray*}
underlies the construction below.
We add markers $\alpha_0$ and $\beta_0$, respectively, to the alphabets
${\bf A}_1$ and ${\bf B}_1$ described above.
Let $B(p)$ be the Bernoulli shift on the alphabet
${\bf A}=\{\alpha_0,\ldots,\alpha_{5}\}=\{\alpha_0\} \cup {\bf A}_1$,
 with associated probability
vector $p=\left({1 \over 2},{1 \over 4},{1 \over 16},{1 \over 16},{1
    \over 16},{1 \over 16}\right)$.  Let $B(q)$ be the Bernoulli shift
on the alphabet ${\bf B}=\{\beta_0,\ldots,\beta_4\}$ with associated
probability vector $q=\left({1 \over 2},{1 \over 8},{1 \over 8},{1
    \over 8},{1 \over 8}\right)$. 

Next we  construct $\Phi$, a finitary isomorphism from $B(p)$ to $B(q)$:

\noindent{\bf Step 0:} If $x_i = \alpha_0$, let $(\Phi(x))_i=\beta_0$; that is,
send markers to markers.

\noindent{\bf Step 1:}  Match the non-marker locations in pairs.  
Suppose that $i$ is paired with $j$. If $\ell_i \neq \ell_j$, we can 
assume that $\ell_i=1$ and $\ell_j=3$ (otherwise reverse the roles).
Remove the bottom bit of $x_j$ and append it below $x_i$, output the
resulting ${\bf B}$ symbols, and {\dof remove from consideration}
both $i$ and $j$.  If $\ell_i=\ell_j$, then
 do not remove $i$ and $j$ from consideration.

\enspace
\noindent For each $n \geq 2$, perform:

\noindent{\bf Step} $n${\bf :}  
The locations which we have not removed from consideration are grouped
in $2^{n-1}$- tuples.  Each such $2^{n-1}$-tuple is either of type 3
(which we define to mean that for every location $i$ within the tuple
$\ell_i=3$), or of type 1. Using the markers,
match the $2^{n-1}$-tuples which have not been removed from
consideration in pairs to form $2^n$-tuples.  If a $2^{n-1}$-tuple
$\xi_3$ of type $3$ is matched with a $2^{n-1}$-tuple $\xi_1$ 
of type $1$, remove the
bottom bit from each $x_i$ in $\xi_3$, and append it to the
corresponding symbol in $\xi_1$.
Finally, output the symbols of ${\bf B}$ thus generated, and remove these
locations from consideration.

\enspace

The coding length for the isomorphism described 
above has essentially the same tails
as Meshalkin's.  To explain this,
observe that the probability $F_k$ that a symbol at the origin is
not coded during the first $k$ pairing stages is approximately
$2^{-k}$ (the approximation is due to parity problems caused by markers.)
After the $k^{{\rm th}}$ pairing stage, only about $1/2^k$ of the symbols
remain uncoded, and these symbols are grouped into $2^k$-tuples.
 Thus heuristically, the event $F_k$ corresponds to an expected coding
distance of order $4^k$. This suggests that
$\py(\n_{\Phi} >t)\approx t^{-{1 \over 2}}$.
Indeed, for this example,  Theorem \ref{strong}
implies that $\E(\n_{\Phi}^{1/2})=\infty$
and the proof of Theorem \ref{one} will show 
that $\E(\n_{\Phi}^{\theta})<\infty$
for all $\theta<1/2$.
\enspace
\enspace

\noindent{\bf An example with $3/4-\epsilon$ moments: heuristics.}
Consider different probability vectors $p$ and $q$,
chosen so that the random walks counting the accumulated information of
non-marker symbols have moment generating functions
\begin{equation} \label{gam3}
\Gamma(z)=\left(\left({1 + z \over 2}\right)^4+
\left({1 - z \over 2}\right)^4\right)z^3
\end{equation}
and
\begin{equation} \label{del3}
\Delta(z)=\left(\left({1 + z \over 2}\right)^4-
\left({1 - z \over 2}\right)^4\right)z^3 \,,
\end{equation}
respectively. Then
\begin{eqnarray*}
\Gamma^2(z)-\Delta^2(z)=
{1 \over 8}\bigg(\Gamma\bigg(z^2\bigg)-\Delta\bigg(z^2\bigg)\bigg).
\end{eqnarray*}
This example is the case $n=2$ of the sequence of examples analyzed in
\S5; see (\ref{gam4}) and (\ref{del4}).

Define a finitary coding $\Phi$ from $B(p)$ to $B(q)$ 
by adapting the recipe above (see \S4 and \S5 for details).
To estimate the tails of $\n_\Phi$, start by observing that
the probability that a symbol is not coded during the first $k$
pairing stages is about $8^{-k}$. At that stage, symbols are grouped
into $2^k$-tuples, and only $1/8^k$ of them remain uncoded,
so heuristically, this event corresponds to an expected
coding distance of order $16^k$. This suggests that
\begin{eqnarray*}
\py(\n_{\Phi} >t)\approx t^{-{3 \over 4}}.
\end{eqnarray*}
Indeed, for this example we will show in \S 5 that
$\E(\n_{\Phi}^\theta)<\infty$ for all $\theta<3/4$. This
is consistent with Theorem \ref{strong}, since the identities
$\Gamma'(1)=\Delta'(1)$ and $\Gamma''(1)=\Delta''(1)$
indicate that $p$ and $q$ 
have the same entropy and the same informational variance.

\section{Ordered measure preserving matchings}

In this section, we define a type of matching which we will employ in
our constructions in \S 5, and derive some useful properties of these
matchings.  Let ${\bf C}=\{\gamma_1,\ldots,\gamma_c\}$ and ${\bf
  D}=\{\delta_1,\ldots,\delta_d\}$ be finite alphabets, and let
$r=(r(\gamma_1),\ldots,r(\gamma_c))$ and
$s=(s(\delta_1),\ldots,s(\delta_d))$ be probability vectors.  Let

\enspace

$\Gamma_k^*=\Gamma(k,{\bf C},r)=\sum_{\bf C}\{r(\gamma_i):r(\gamma_i)=2^{-k}\}$

\enspace

\noindent and

\enspace

$\Delta_k^*=\Delta(k,{\bf D}, s)=\sum_{\bf D}\{s(\delta_j):
s(\delta_j)=2^{-k}\}$.

\enspace

Define an order relation $\prec$ on ${\bf C}$ such that $\gamma_1
\prec \cdots \prec \gamma_c$ and an order relation $\prec$ on ${\bf
  D}$ such that $\delta_1 \prec \cdots \prec \delta_d$. 
Endow ${\bf C} \times {\bf C}$ with the lexicographic ordering, i.e., define
$\gamma_i\gamma_j \prec \gamma_{m}\gamma_{n}$ if $\gamma_i \prec
\gamma_{m}$ or if $\gamma_i=\gamma_{m}$ and $\gamma_j \prec
\gamma_{n}$.  Similarly, endow ${\bf D} \times {\bf D}$
with the  lexicographic ordering $\prec$.

Let $r(\gamma_i\gamma_j)=r(\gamma_i)r(\gamma_j)$ and
$s(\delta_i\delta_j)=s(\delta_i)s(\delta_j)$.  We define the {\dof
  maximal ordered measure preserving matching} ({\dof mompm})
$\psi=\psi_{({\bf C},{\bf D},r,s)}$ from ${\bf C} \times {\bf C}$ to
${\bf D} \times {\bf D}$ given $(r,s)$ as follows:

For all $t \in {\mathbb R}$, write the ordered set $\{x \in {\bf C} \times
{\bf C}:r(x)=t\}$ in increasing order as $\{x_t(i):1 \leq i \leq
\ell_t\}$, and similarly, write the ordered set $\{y \in {\bf D}
\times {\bf D}:s(y)=t\}$ in increasing order as $\{y_t(i):1 \leq i
\leq m_t\}$, assuming these sets are non-empty.  Define
$\psi(x_t(i))=y_t(i)$ for $1\leq i \leq \min\{\ell_t,m_t\}$.

Let $E=E({\bf C},{\bf D},r,s)$ be the set in ${\bf C} \times {\bf C}$
where $\psi$ is defined.  Let $F=F({\bf C},{\bf D},r,s)=\psi(E)$.  Let
$G=G({\bf C},{\bf D},r,s)={\bf C} \times {\bf C}-E$.  Let $H=H({\bf
  C},{\bf D},r,s)={\bf D} \times {\bf D}-F$.

Let $\widetilde{r}={\widetilde{r}}_{({\bf C},{\bf D},r,s)}=
\left(\widetilde{r}(x):x\in G\right)$, where
$\widetilde{r}(x)={r(x)\over \sum_{\widetilde{x} \in
    G}r(\widetilde{x})}$ and $\widetilde{s}={\widetilde{s}}_{({\bf
    C},{\bf D},r,s)}=\left(\widetilde{s}(y):y \in H\right)$, where
$\widetilde{s}(y)={s(y)\over \sum_{\widetilde{y} \in
    H}s(\widetilde{y})}$ be the {\dof probability vectors induced} by
$(r,s)$ on $G$ and $H$.

Let

\enspace

$\Upsilon_k^*= \Upsilon(k,{\bf C},r)=\sum_{x \in {\bf C} \times {\bf
    C}}\{r(x):r(x)=2^{-k}\}$,

\enspace

$\Omega_k^*= \Omega(k,{\bf D},s)=\sum_{y \in {\bf D}\times{\bf
    D}}\{s(y):s(y)=2^{-k}\}$,

\enspace

$\Lambda_k^*=\Lambda(k,{\bf C},{\bf D},r,s)=\sum_{x \in
  G}\{r(x):r(x)=2^{-k}\}$,

\enspace

\noindent and let

\enspace

$\Xi_k^*=\Xi(k,{\bf C}, {\bf D},r,s)=\sum_{y \in H}\{s(y):s(y)=2^{-k}\}$.

\enspace
\enspace

We say that $\psi$ {\dof reduces mass} by a factor of $t$ if
$\sum_{k=0}^{\infty}\Lambda_k^*=t$.

Let 
\begin{eqnarray}
\label{gamma}
\Gamma(z)=\Gamma({\bf C},{\bf D},r,s,z)=\sum_{k=0}^\infty\Gamma_k^*z^k. 
\end{eqnarray}

Define $\Delta(z)$, $\Upsilon(z)$,
$\Omega(z)$, $\Lambda(z)$, and $\Xi(z)$ analogously.  Then $\Upsilon(z)=\Gamma^2(z)$ and $\Omega(z)=\Delta^2(z)$.
Also, $\Lambda(z)-\Xi(z)=\Upsilon(z)-\Omega(z)$.

\begin{lemma}
\label{ladder}
Suppose $\Gamma^2(z)-\Delta^2(z)=t\Gamma(z^2)-t\Delta(z^2)$.  Then:
\begin{mylist}
\item $\Lambda(z)=t\Gamma(z^2)$ and $\Xi(z)=t\Delta(z^2)$ 
\item $\psi$ reduces mass by a factor of $t$.
\end{mylist}
\end{lemma}
\begin{pf}
\begin{mylist}
\item $\Lambda(z)-\Xi(z)=\Upsilon(z)-\Omega(z)=\Gamma^2(z)-\Delta^2(z)=t\Gamma(z^2)-t\Delta(z^2)$,

\noindent hence 
\begin{eqnarray}
\Lambda(z)=t\Gamma(z^2) \label{lamb}
\end{eqnarray}
and
\begin{eqnarray}
\Xi(z)=t\Delta(z^2). \label{xi}
\end{eqnarray}
\item
By (\ref{lamb}),
\begin{eqnarray*}
\sum_{k=0}^\infty\Lambda_k^*=t\sum_{k=0}^{\infty}\Gamma_k^*=t.
\end{eqnarray*}
\end{mylist}
\end{pf}

Let $C_1={\bf C}$, let $D_1={\bf D}$, let $r_1=r$, and let $s_1=s$.
Inductively, let $C_{i+1}=G(C_i,D_i,r_i,s_i)$, let
$D_{i+1}=H(C_i,D_i,r_i,s_i)$, let
$r_{i+1}=\widetilde{r}_{(C_i,D_i,r_i,s_i)}$, and let
$s_{i+1}=\widetilde{s}_{(C_i,D_i,r_i,s_i)}$.  Let
$\psi_i=\psi{(C_i,D_i,r_i,s_i)}$.  Note that $\psi_i$ matches
$2^i$-tuples to $2^i$-tuples.  We call $\{\psi_i\}_{i \ge 1}$ the
{\dof sequence of mompm's associated to} $({\bf C},{\bf D},r,s)$.  Let
$\Gamma_i(z)=\Gamma(C_i,D_i,r_i,s_i,z)$.  In particular,
$\Gamma_1(z)=\Gamma(z)$ as defined in equation (\ref{gamma}).  Define
$\Delta_i(z)$, $\Upsilon_i(z)$, $\Omega_i(z)$, $\Lambda_i(z)$, and
$\Xi_i(z)$ analogously.

Inductive application of Lemma \ref{ladder} gives:

\begin{cor}
\label{matchup}
Suppose $\Gamma^2(z)-\Delta^2(z)=t\Gamma_1(z^2)-t\Delta_1(z^2)$.
\begin{mylist}
\item  If $i \in {\mathbb Z}_+$, then $\Gamma_i^2(z)-\Delta_i^2(z)=t\Gamma_i(z^2)-t\Delta_i(z^2)$
\item  If $i \in {\mathbb Z}_+$, then $\psi_i$ reduces mass by a factor of $t$.
\end{mylist}
\end{cor}

\section{A class of codes with finite moments}

Finally, we construct a class of examples to prove Theorem \ref{one}.

Fix $n \in {\mathbb Z}_+$.

Let $p_0=q_0={1 \over 2}$.  
Construct $p=({1 \over 2},p_1,\ldots,p_{a-1})$ such that
for each integer $m \in[0,n]$, 
exactly $2^{2m}{2n \choose 2m}$ of the $p_i$ take the value $2^{-2m-2n}$.
Thus if for $i \geq 1$, we 
denote $r(\alpha_i)=2p_i$, then  
\begin{eqnarray}
\label{odd}
\Gamma_{2m+2n-1}^*= \sum_i\{r(\alpha_i):r(\alpha_i)=
2^{-(2m+2n-1)}\}={1 \over 2^{2n-1}}{2n \choose 2m}
\end{eqnarray}
for all $m \in {\mathbb Z}$ such that $0 \leq m \leq n$. Define 
$\Gamma_k^* =0$ for all other $k$.

Similarly, for $j \geq 1$, denote $s(\beta_j)=2q_j$ and 
construct $q=({1 \over 2}, q_1,\ldots,q_{b-1})$ such that
\begin{eqnarray}
\label{even}
\Delta_{2m+2n}^*= \sum_j\{s(\beta_j):s(\beta_j)=2^{-(2m+2n)}\}={1 \over 2^{2n-1}}{2n \choose 2m+1}
\end{eqnarray}
for all $m \in {\mathbb Z}$ such that $0 \leq m \leq n-1$, and define 
$\Delta_k^* =  0$ for other $k$.

Let $B(p)$ be the Bernoulli shift with probability vector $p$ on the alphabet ${\bf A}=\{\alpha_0,\ldots,\alpha_{a-1}\}$.
Let $B(q)$ be the Bernoulli shift
 with probability vector $q$ on the alphabet 
${\bf B}=\{\beta_0,\ldots,\beta_{b-1}\}$.

Let ${\bf C}=\{\alpha_1,\ldots,\alpha_{a-1}\}$ and let ${\bf
  D}=\{\beta_1, \ldots, \beta_{b-1}\}$. Consider the probability
  vectors
$r=(r(\alpha_i):1 \leq i
\leq a-1)$ and let $s=(s(\beta_j):1 \leq j \leq b-1)$. Relative to
these, define all other terms as in \S 4.

\begin{lemma}
If $i \in {\mathbb Z}_+$, then $\psi_i$ reduces mass by a factor of ${1 \over 2^{2n-1}}$.
\end{lemma}
\begin{pf}  
Recall that $\Gamma(z)=\sum_{k=0}^\infty\Gamma_k^*z^k$ and $\Delta(z)=\sum_{k=0}^\infty\Delta_k^*z^k$.  By the binomial theorem and equations (\ref{odd}) and
(\ref{even}), we find that
\begin{eqnarray*}
\Gamma^2(z)-\Delta^2(z) &=& (\Gamma(z)-\Delta(z))(\Gamma(z)+\Delta(z))\\
&=& 2\Bigl({1-z \over 2}\Bigr)^{2n}z^{2n-1}2\Bigl({1+z \over 2}\Bigr)^{2n}z^{2n-1}\\
&=& {1 \over 2^{2n-1}}2\Bigl({1-z^2 \over 2}\Bigr)^{2n}z^{4n-2}\\
&=& {1 \over 2^{2n-1}}\Bigl( \Gamma\Bigl(z^2\Bigr)-\Delta\Bigl(z^2\Bigr)\Bigr),
\end{eqnarray*}
so the desired result holds by Corollary \ref{matchup}.
\end{pf}

\noindent {\bf Example}

When $n=2$, we may let ${\bf A}= \{ \alpha_0, \ldots , \alpha_{41}\}$; ${\bf B}= \{ \beta_0, \ldots , \beta_{40}\}$; 
$p=(p_0, \ldots , p_{41})$ such that $p_0 = 2^{-1}$, $p_1 = 2^{-4}$, $p_2= \cdots = p_{25} = 2^{-6}$, and 
$p_{26}= \cdots = p_{41} = 2^{-8}$; and
$q=(q_0, \ldots , q_{40})$ such that $q_0 = 2^{-1}$, $q_1= \cdots = q_{8} = 2^{-5}$, and $q_{9}= \cdots = q_{40} = 2^{-7}$.

Taking logarithms to base $2$, $h(p)=h(q)={7 \over 2}$ and $\sigma^2_p
= \sigma^2_q={27 \over 4}$, hence Theorem \ref{strong} does not apply.
These vectors correspond to the generating functions
in (\ref{gam3}) and (\ref{del3}). We find that

\begin{equation} \label{gam4}
\Bigl(\Gamma_3^*, \Gamma_5^*, \Gamma_7^*\Bigr)= 
\Bigl({1 \over 8}, {3 \over 4},{1 \over 8}\Bigr)
\end{equation}

\begin{equation} \label{del4}
\Bigl(\Delta_4^*, \Delta_6^*\Bigr)=\Bigl({1 \over 2},{1 \over 2}\Bigr)
\end{equation}

\begin{equation}
\Bigl(\Upsilon_6^*,\Upsilon_8^*,\Upsilon_{10}^*,
\Upsilon_{12}^*,\Upsilon_{14}^*\Bigr)=\Bigl({1
  \over 64},{3 \over 16},{19 \over 32},{3 \over 16},{1 \over
  64}\Bigr)
\end{equation}

\begin{equation}
\Bigl(\Omega_8^*, \Omega_{10}^*, \Omega_{12}^*\Bigr)=
\Bigl({1 \over 4},{1 \over 2},{1 \over 4}\Bigr)
\end{equation}

\begin{equation}
\Bigl(\Lambda_{6}^*,\Lambda_{10}^*,\Lambda_{12}^*\Bigr)=
\Bigl({1 \over 64},{3 \over 32},{1 \over 64}\Bigr)
\end{equation}

\begin{equation}
\Bigl(\Xi_8^*,\Xi_{12}^*\Bigr)=\Bigl({1 \over 16},{1 \over 16}\Bigr).
\end{equation}

\noindent {\bf Definition of $\Phi$}

For $x =(x_k)_{k \in {\mathbb Z}} \in X = {\bf A}^{\mathbb Z}$, define a {\dof $j$-marker} as a run of at least $2nj$ consecutive $\alpha_0$ symbols.  Define a 
{\dof $j$-gap} as the location of the non-$\alpha_0$ symbols between neighboring $j$-markers.

Let $G_{j,0}=\{g(j,0,1),\ldots,g_(j,0,\ell_{j,0})\}$ be the
ordered elements (from left to right) of the $j$-gap containing
$\min\{i \geq 0:x_i \neq \alpha_0\}$.  More generally, let
$G_{j,i}=\{g(j,i,1),\ldots,g(j,i,\ell_{j,i})\}$ be the ordered
elements of the $i^{th}$ $j$-gap to the right of $G_{j,0}$ (to the
left if $i <0$).

\enspace

{\dof Step 0:}  If $x_i = \alpha_0$, let $(\Phi(x))_i = \beta_0$.

{\dof Step 1:}  Within each 1-gap, match the elements in pairs, starting from the left ($g(1,i,1)$ with $g(1,i,2)$, $g(1,i,3)$ with $g(1,i,4)$, etc.).
All the elements will be paired except possibly $g(1,i,\ell_{1,i})$.  

If $\psi_1(x_{g(1,i,2k+1)}x_{g(1,i,2k+2)})$ is defined, then let

\enspace

$(\Phi(x))_{g(1,i,2k+1)}(\Phi(x))_{g(1,i,2k+2)}=
\psi_1(x_{g(1,i,2k+1)}x_{g(1,i,2k+2)})$,

\enspace

\noindent and {\dof remove from consideration} 
$g(1,i,2k+1)$ and $g(1,i,2k+2)$.  

Starting from the left, match the pairs which have not been removed from consideration into quartets.  
If $\psi_2$ of the symbols at the position of a quartet
is defined, output the result in the position of the quartet and remove the elements of the quartet from consideration.  

Iterate, matching
$2^{k-1}$-tuples which have not been removed from consideration into $2^k$-tuples and applying $\psi_k$, until $2^k > \ell_{1,i}$.

\enspace

For each $j \geq 2$, do the following:

\enspace

{\dof Step} $j${\bf :}  Within each $j$-gap, starting from the left, match into pairs any elements in $G_{j,i}$ which were not paired in any of the previous steps,
and apply $\psi_1$ as in Step 1.  

Match into quartets any previously unmatched pairs (including the
pairs just created) which have not been removed from consideration,
and apply $\psi_2$, etc., iterating until $2^k > \ell_{j,i}$.

\enspace
\enspace

When $n=1$, this is the code described in \S 3.
The next two lemmas are needed as preparation for
bounding the tails of $\n_\Phi$.

\begin{lemma}\label{smallell}
If $f(x)=\ell_{j,0}^{-1}\ind_{[x_0 \ne \alpha_0]} $, then 
$\E(f) = 2^{-2nj-1}$.
\end{lemma}
\begin{pf}
The sum $\sum_{m=1}^M f(T^m x)$
differs from the number of $j$-gaps in $[1,M]$ by at most $2$.
Counting $j$-gaps in $[1,M]$ is equivalent to counting runs of $2nj$ marker
symbols followed by a non-marker symbol; such strings have asymptotic
frequency $2^{-2nj-1}$. Taking the limit of 
$\frac{1}{M}\sum_{m=1}^M f(T^m x)$
as $M \to \infty$,
the ergodic theorem yields the assertion.
\end{pf}

\begin{lemma}
\label{bigL}
For $j \ge 1$ 
Let $L_{j,i}=2nj+g(j,i,\ell_{j,i})-g(j,i,1)$ denote the ``span''
of the $i^{th}$ $j$-gap.
  If $\theta < 1$, then 
$$
\E\Bigl((L_{j,0}-L_{j-1,0})^\theta \mid x_0 \neq \alpha_0\Bigr) 
\le 2^{(2+2nj)\theta}\,.
$$

\end{lemma}
\begin{pf}
The expected distance between the beginnings of successive $j$-gaps
is $2^{1+2nj}$ by Kac's Theorem (see \cite{pete}, p.\ 46), 
whence 
$$
\E\Bigl(L_{j,0}-L_{j-1,0} \mid x_0 \neq \alpha_0\Bigr) \le 
2^{2+2nj} \,.
$$
The assertion of the lemma follows by Jensen's inequality.
\end{pf}

\begin{lemma}
\label{forward}
If $\theta < 1- {1 \over 2n}$, then $\E(\n_{\Phi}(x))^\theta < \infty$.
\end{lemma}
\begin{pf}
Recall $L_{j,i}$ from the previous lemma and define $L_{0,i}=0$.
If Step $j$ determines $(\Phi(x))_0$, then $\n_\Phi(x) \leq L_{j,0}$.  
Let $A_j$ be the event that Steps 1 to $j$ do not determine $(\Phi(x))_0$. 

Let $B_j$ be the event that the $0^{th}$ coordinate is 
matched at least $j$ times by the end of Step $j$, but $(\Phi(x))_0$ has not
yet been determined.  
Let $C_j$ be the event that at the end of Step $j$, 
the $0^{th}$ coordinate has been matched at most $j-1$ times
(so it is not part of a $2^j$-tuple).
Clearly, for each $j \geq 1$,
\begin{equation} \label{abc}
\py(A_j) \le \py(B_j)+\py(C_j) \,.
\end{equation}

Every time an undetermined coordinate is matched, the probability 
that it remains undetermined is 
${1 \over 2^{2n-1}}$, whence
\begin{equation} \label{eq:b}
\py(B_j) \le \Bigl({1 \over 2^{2n-1}}\Bigr)^j  \,.
\end{equation}

Since, for all $k$ and $j$, at most one $2^k$-tuple in $G_{j,0}$
is unmatched at the end of Step $j$, it follows that
$$
\py(C_j) \leq
\E\Bigl({\sum_{k=0}^{j-1}2^k \over \ell_{j,0}} \mid x_0 \neq
\alpha_0\Bigr) \leq   2^j\Bigl({1 \over 2^{2n}}\Bigr)^j 
=\Bigl({1 \over 2^{2n-1}}\Bigr)^j
$$
by Lemma \ref{smallell}. Thus
\begin{eqnarray*}
\py(A_j) &\leq& 
 2\Bigl({1 \over 2^{2n-1}}\Bigr)^j.
\end{eqnarray*}

Therefore
\begin{eqnarray*}
\E(N_\Phi(x))^\theta &\leq& \sum_{j=1}^\infty \py(A_{j-1})
\E\left(L_{j,0}^\theta-L_{j-1,0}^\theta \mid A_{j-1}\right)\\
&\leq& \sum_{j=1}^\infty 2\left({1 \over 2^{2n-1}}\right)^j 
\E\left((L_{j,0}-L_{j-1,0})^\theta \mid A_{j-1}\right)\\
\end{eqnarray*}

Conditional on the event that $x_0 \neq \alpha_0$, the random variable
$(L_{j,0}-L_{j-1,0})$ is independent of the event $A_{j-1}$, hence
by Lemma \ref{bigL},
\begin{eqnarray*}
\E(N_\Phi(x))^\theta &\leq& \sum_{j=1}^\infty 
2\Bigl({1 \over 2^{2n-1}}\Bigr)^j \E\Bigl((L_{j,0}-L_{j-1,0})^\theta 
\mid x_0 \neq \alpha_0\Bigr)\\
&\leq& \sum_{j=1}^\infty 2 \Bigl({1 \over 2^{2n-1}}\Bigr)^j 
4(2^{2n})^{j\theta} 
= 8 \sum_{j=1}^\infty \Bigl({1 \over 2^{2n-1}} \Bigr) ^{j\Big
  (1-{2n\theta \over 2n-1}\Bigr ) } < \infty.\\
\end{eqnarray*} 
\end{pf}

A similar argument gives:

\begin{lemma}
\label{backward}
If $\theta < 1- {1 \over 2n}$, then $\E_Q(\n_{\Phi^{-1}}(x))^\theta < \infty$.
\end{lemma}

\begin{pfof}{Theorem \ref{one}}
  By Lemmas \ref{forward} and \ref{backward}, it only remains to
  verify that $\Phi$ is an isomorphism.  Since $\Phi$ is finitary, it
  gives an a.e. defined map from $B(p)$ to $B(q)$. 
As our   definition of $(\Phi(x))_i$ depends only on the position of $i$
  within its $j$-blocks, $\Phi$ is translation invariant.
Since each
  $\psi_k$ is a one-to-one measure preserving matching from
  previously uncoded sequences to previously uncoded sequences, it
  follows that $\Phi$ is measure preserving and invertible.  
More precisely, for $\P$-a.e.\ $x \in X$ and any $n \ge 1$,
all the symbols in the string $(x_{-n}, \ldots, x_n)$
get coded within a finite distance. This means that the cylinder set
$\Bigl\{\widetilde{x} \in X : 
(\widetilde{x}_{-n}, \ldots, \widetilde{x}_n)  =  
(x_{-n}, \ldots, x_n) \Bigr\}$ 
is partitioned into countably many cylinder sets ${\mathcal C}_j$
(and a set of measure zero); each ${\mathcal C}_j$
is mapped, using one of our matchings $\psi_{k(j)}$, 
to a cylinder set $\Phi({\mathcal C}_j)$ 
in $Y$ with $\Q(\Phi({\mathcal C}_j))=\P({\mathcal C}_j)$. 
This completes the proof.

\end{pfof}

\section{Extension to ergodic Markov chains}

Let ${\bf A}= \{\alpha_0, \ldots ,\alpha_{a-1}\}$ be a finite alphabet
and let $p=(p(\alpha_i,\alpha_j))_{0 \leq i,j\leq a-1}$ be an irreducible
 stochastic matrix.  The associated Markov chain
$M(p)$ is ergodic and has a (strictly positive) 
unique stationary distribution
$\widetilde{p}=(\widetilde{p}(\alpha_0),\ldots,\widetilde{p}(\alpha_{a-1}))$.
  Similarly, let ${\bf B}=
\{\beta_0, \ldots ,\beta_{b-1}\}$ be a finite alphabet and let
$q=(q(\beta_i,\beta_j))_{0 \leq i,j \leq b-1}$ be a stochastic matrix
such that $M(q)$ is ergodic with unique stationary distribution
$\widetilde{q}=(\widetilde{q}(\beta_0),\ldots,\widetilde{q}(\beta_{b-1}))$.
The Markov chain $M(p)$ has entropy 
$$
h(p)=\sum_{0 \leq i,j \leq
a-1}-\widetilde{p}(\alpha_i)p(\alpha_i,\alpha_j)\log p(\alpha_i,\alpha_j) \,.
$$
We will assume that $h(p) = h(q)$. 
Let $\varphi$ be a finitary homomorphism from $M(p)$ to $M(q)$.  For
$x=(x_k)_{k \in {\mathbb Z}}\in {\bf A}^{\mathbb Z}$, let $X_i(x)=-\log
(p(x_{i-1},x_i))-h(p)$.  Similarly, if $\varphi(x)=y=(y_k)_{k \in
{\mathbb Z}}\in {\bf B}^{\mathbb Z}$, let $Y_i(x)=-\log(q(y_{i-1},y_i))-h(q)$.
Let $S_{m,n}=\sum_{i=m+1}^nX_{i}$ and let
$R_{m,n}=\sum_{i=m+1}^nY_{i}$.  Since $\varphi$ is measure preserving,
it follows that ${\bf E}(X_i)={\bf E}(Y_i)=0$.  Let 
$$
\lambda_p=\max_{0
\leq i,j \leq a-1}\{-\log(p(\alpha_i,\alpha_j)):p(\alpha_i,\alpha_j)
\neq 0\} \, ,
$$
and let $\gamma_p=\max_{0 \leq i\leq
a-1}\{-\log(\widetilde{p}(\alpha_i))\}$.

The following central limit theorem can be found, e.g., in
 \cite{durrett}, p.\ 422 under an additional aperiodicity assumption,
and in \cite{GL} in much greater generality.
For the reader's convenience, we include a brief proof.
\begin{lemma}
If $M(p)$ is an ergodic Markov chain on a finite alphabet, then there
exists a constant $\sigma_p \geq 0$ depending only on $p$ such that
${S_{0,n} \over \sqrt{n}}\Rightarrow \chi\sigma_p$ in law,
where $\chi$ denotes a standard normal variable.
\end{lemma}
We define $\sigma^2_p$ to be the {\dof 
asymptotic informational variance} of $p$.

\begin{pf}
For any $x \in {\bf A}^{\mathbb Z}$, let $T_0=\min\{t>0:x_t=\alpha_0\}$.
Inductively, for $i \geq 0$, let $T_{i+1}=\min\{t>T_i:x_t=\alpha_0\}$.
The increments $T_{i}-T_{i-1}$ are i.i.d.\ and have exponential tails.
The partial sums $\{S_{T_{i-1},T_i}\}_{i \geq 1}$ are also i.i.d.
Let $d_p={\bf E}(T_1-T_0)>0$.
 By an application of the
ergodic theorem and the law of large numbers, 
${\bf E}(S_{T_0,T_1})=0$.  
Since $|X_i|\leq \lambda_p$, it follows that
$S^2_{T_0,T_1}\leq (T_1-T_0)^2\lambda_p^2$, whence 
${\bf E}(S^2_{T_0,T_1})=c_p^2<\infty$.
Let $N_n=\min\{m>0:T_m \geq n\}$.  Since ${N_n \over n/d_p}
\rightarrow 1$ in probability, the random index central
limit theorem (see \cite{durrett}, p. 116) states that
\begin{equation} \label{sub}
{S_{T_0,T_{N_n}} \over \sqrt{n}}\Rightarrow {c_p \chi \over
\sqrt{d_p}}.
\end{equation}
Define  $\sigma^2_p={c^2_p \over d_p}$.
Since 
$
{\bf E}(T_{N_n}-n)\leq \max_{0 \leq i \leq a-1}{\bf E}(T_0\mid x_0=\alpha_i)
$
 for all $n \in {\mathbb Z}_+$, 
it follows that
\begin{eqnarray*}
{\bf E}(|S_{0,n}-S_{T_0,T_{N_n}}|)\leq \lambda_p{\bf E}(T_0)+\lambda_p
\max_{0 \leq i \leq a-1}{\bf E}(T_0\mid x_0=\alpha_i).
\end{eqnarray*}
In conjunction with (\ref{sub}), this gives
${S_{0,n} \over \sqrt{n}} \Rightarrow \chi\sigma_p $.
\end{pf}

Let $J_n = \{i \in \{1,\ldots, n\}:n_{\varphi}(T^ix) >
\min\{i,n+1-i\}$ or $n_{\varphi}(T^{i-1}x) > \min\{i-1,n+2-i\}\}$.
Let $I_n=\{1, \ldots,n\}-J_n$.

As in \S 2, we deduce from the CLT and uniform integrability:
\begin{lemma}
If $\sigma^2_p \neq \sigma^2_q$, then $\liminf_{n \to \infty}
\frac{1}{\sqrt n}\E (R_{0,n}-S_{0,n})^+ \ge
\frac{|\sig_q-\sig_p|}{\sqrt{2\pi}}.$
\end{lemma}
The proofs of Lemma \ref{short} and Theorem \ref{strong} 
adapt to prove the following.
\begin{lemma}
Suppose $M(p)$ and $M(q)$ are ergodic Markov chains and $\varphi$ is a
finitary homomorphism from $M(p)$ to $M(q)$, Then for all $n$,
\begin{eqnarray*}
{\bf E}(R_{0,n}-S_{0,n})^+\leq \gamma_p+4\lambda_q({\bf
E}(N_\varphi\wedge n)+1)
\end{eqnarray*}
\end{lemma}
\begin{thm}
Let $M(p)$ and $M(q)$ be ergodic Markov chains such that $h(p)=h(q)$
and $\sigma^2_p \neq \sigma^2_q$. If $\varphi$ is a finitary
homomorphism from $M(p)$ to $M(q)$, then ${\bf
E}\left(\sqrt{N_\varphi(x)}\right)= \infty$.  More precisely, $
\liminf_{n \rightarrow \infty} {1\over \sqrt{n}} {\bf E} (
N_\varphi(x) \wedge n) \geq c_{p,q}>0$.
\end{thm}

\section{Higher moments: a problem}
Theorem \ref{strong} and our constructions in \S 5
suggest the following:
\enspace

\noindent{\bf Question.}
Let $p$ and $q$ be probability vectors with
$h(p)=h(q)$. Fix an integer $ k > 2$. Suppose that $\varphi$ is a 
finitary homomorphism  from $B(p)$ to $B(q)$,
that satisfies $\E\Bigl(\n_\varphi^{1-1/k}\Bigr)<\infty$.
Does it follow that
$$
\sum_i p_i (\log p_i)^k = \sum_j q_j (\log q_j)^k \; \: ?
$$
\section*{Acknowledgment}
We thank Alexander Holroyd for helpful discussions,
Serban Nacu for comments and Ben-Zion Rubshtein 
and Jeff Steif for references.

\end{document}